\documentclass[12pt]{amsart}
\newtheorem{theorem}{Theorem}

\begin{document}

\title{A converse of the Gauss--Lucas theorem}
\author{Nikolai Nikolov and Blagovest Sendov}

\address{Institute of Mathematics and Informatics\\Bulgarian Academy of
Sciences\\1113 Sofia, Bulgaria}\email{nik@math.bas.bg}

\address{Institute of Information and Communication Technologies\\Bulgarian Academy of Sciences\\
1113 Sofia, Bulgaria}\email{sendov2003@yahoo.com}

\begin{abstract}
 All linear operators $L:\mathcal{C}[z]\to
\mathcal{C}[z]$ which decrease the diameter of the zero set of any
$P\in\mathcal{C}[z]$ are found.
\end{abstract}

\maketitle

Let $Q\in\mathcal{C}[z]$ be a complex polynomial with zero set $Z(Q)$. One
of the most important and useful theorems in the Geometry of
polynomials is the  Gauss-Lucas theorem (see, for example,
\cite{MM1}, \cite{RS}).

\begin{theorem}\label{GL} {\bf(Gauss-Lucas theorem)} If $P\in\mathcal{C}[z]$
is not a constant, then the convex hull of $Z(P)$ contains $Z(P')$.
\end{theorem}

This classical theorem has very important physical interpretation
in the potential theory and various proofs, some of them rather
simple. It may be considered as a variant of the Rolle's theorem
for complex polynomials.

In \cite{GS}, a converse of Theorem \ref{GL} is proved, which we stated as follows

\begin{theorem}\label{GS}
Let $L:\mathcal{C}[z]\to\mathcal{C}[z]$ be a linear operator such that if $L[P]\neq
0,$ then the convex hull of $Z(P)$ contains $Z(L[P]).$ Then either
$L$ is a linear functional, or there exist
$c\in\mathcal{C}_\ast=\mathcal{C}\setminus\{0\}$ and $k\ge 0$ such that
$L[P]=cP^{(k)}$ (the $k$th derivative).
\end{theorem}

Set $\mathcal{P}_n=\{P\in\mathcal{C}[z]:\deg P\leq n\}.$  Theorem \ref{GS} fully
characterized the linear operators $L:\mathcal{C}[z]\to\mathcal{C}[z]$ with the
property that the convex hull of $Z(P)$ contains $Z(L[P]).$ A big
chapter in the Geometry of polynomials is devoted to linear
operators $L:\mathcal{C}[z]\to\mathcal{C}[z]$ with given properties, defined by
specific geometric relation between $Z(L[P])$ and $Z(P)$. For a
recent development of this topic see, for example, \cite{BB} and
\cite{CC}.

The proof of Theorem \ref{GS} uses only the linearity of $L$ and the
fact that the convex
hull of a singleton is the singleton itself, which implies that
$L[(z+\alpha)^n]=a_n(\alpha)(z+c_n(\alpha))^{k_n(\alpha)}$.

Modifying the approach from \cite{GS}, we shall prove a more
general result in the same spirit.  For $0\le k\le n-2$ denote by
$d_{n,k}\in(0,1]$ the smallest number such that
diam$Z(P^{(k)})\le d_{n,k}$ diam$Z(P)$ for any $P\i\mathcal{P}_n\setminus\mathcal{P}_k$
(where diam$\emptyset=0$). It is not difficult to see that
$d_{n,k}=1$ if and only if $2k\le n-2.$ Note also that $\displaystyle
d_{3,1}=\frac{2}{3}.$

\begin{theorem}\label{NS} A linear operator $L:\mathcal{P}_n\to\mathcal{C}[z]$
($n\ge 1$) shares the property that diam$Z(L[P])\le$ diam$Z(P),$ when both $P$ and $L[P]$ are not
constants,\footnote{Assuming diam$\emptyset=0,$ one may
replace this property by diam$Z(L[P])\le$ diam$Z(P)$ if
$L[P]\neq 0.$ Then the first argument in the proof is
superfluous.} if and only if $L$ has one of the forms:

1) $L=zl_1+l_2,$ where $l_1$ and $l_2$ are linear functionals;

2) $L=(z-c)^ml_3,$ where $m\ge 2,$ $c\in\mathcal{C}$ and $l_3\neq 0$ is a
linear functional;

3) there exist $c\in\mathcal{C}_\ast,$ $0\le k\le n-2$ and $\mathcal{L}\in\mathcal{P}_1$ with
$|\mathcal{L}'|\ge d_{n,k}$ such that $L[P]=c(P\circ\mathcal{L})^{(k)}.$\footnote{
$P\circ L(z):=P(L(z)).$ If $n=\infty$ (i.e. $L:\mathcal{C}[z]\to\mathcal{C}[z]$),
then $|\mathcal{L}'|\ge 1.$}
\end{theorem}

{\bf Proof.} If $L$ has one of the three forms, then $L$
obviously satisfies the given condition.

To prove the converse, note that for any $\alpha\in\mathcal{C}_\ast,$ $L[\alpha
z+1]=\alpha L[z]+L[1]$ either vanishes, or has no distinct zeros.
Letting $\alpha\to 0,$ the same follows for $L[1].$

Assume now that $L$ is not a linear functional, i.e. $\deg
L[z^r]\ge 1$ for some $r,$ $0\le r\le n.$ Let $r$ be minimal,
$r\le s\le n,$ and $\alpha\in\mathcal{C}.$ Then
$$(a_s(\alpha)z+b_s(\alpha))^{q_s(\alpha)}=L[(z+\alpha)^s]=
\sum_{j=0}^s\binom{s}{j}\alpha^j L[z^{s-j}].$$ Considering the last
sum as a polynomial of $z,$ we see that there exists a finite set
$A_s\subset\mathcal{C}$ such that the degree of this polynomial is positive
and does not depend on $\alpha\not\in A_s.$ After translation, we may
assume that $\displaystyle 0\not\in U=\cup_{j=r}^n A_j.$ Let $\alpha\not\in U.$
It follows that $q_s(\alpha)=q_s(0)=q_s\ge 1$ and $a_s(\alpha)\neq 0.$
Then
$$(a_r(\alpha)z+b_r(\alpha))^{q_r}=
(a_r(0)z+b_r(0))^{q_r}+Q_r(\alpha).$$ In particular,
$a_r(\alpha)=a_r(0)=a_r.$

Let $q_r\ge 2.$ Then $b_r(\alpha)=b_r(0)=b_r$ and $Q_r(\alpha)=0,$ i.e.
$L[1]=\dots=L[z^{r-1}]=0.$ If $r=n,$ then $L$ has the second form.
Otherwise,
$$(a_{r+1}(\alpha)z+b_{r+1}(\alpha))^{q_{r+1}}=
(a_{r+1}(0)z+b_{r+1}(0))^{q_{r+1}}+\alpha(r+1)(a_rz+b_r)^{q_r}.$$
Note that $q_{r+1}\ge q_r.$ Set $z_0=-b_r/a_r.$ If
$(a_{r+1}(0)z+b_{r+1}(0))^{q_{r+1}}\neq 0$ and $\alpha\neq 0,$ then
$z_0$ is a zero of
$$1-\left(\frac{a_{r+1}(\alpha)z+b_{r+1}(\alpha)}
{a_{r+1}(0)z+b_{r+1}(0)}\right)^{q_{r+1}}$$
with multiplicity $q_r.$ The last function of $z$ either vanishes
or has only simple zeros which contradicts that $q_r\ge 2.$ Hence,
$a_{r+1}(0)z_0+b_{r+1}(0)=0$ and $z_0$ is zero of
$(a_{r+1}(\alpha)z+b_{r+1}(\alpha))^{q_{r+1}}$ with multiplicity $q_r.$
It follows that $q_{r+1}=q_r$ and $L[z^{r+1}]=c_{r+1}L[z^r].$
Repeating the same arguments, we conclude that $L$ has the second
form with $c=z_0$ and $m=q_r.$

It remains to deal with the case $q_r=1.$ Assume that $L$ is not
of the first form. Let $q_r=\dots=q_t=1$ and $q_{t+1}\ge 2.$ Then
$$(a_{t+1}(\alpha)z+b_{t+1}(\alpha))^{q_{t+1}}=(a_{t+1}(0)z+b_{t+1}(0))^{q_{t+1}}
+\sum_{j=1}^{t+1}\binom{t+1}{j}\alpha^jL[z^{t+1-j}].$$

If $q_{t+1}\ge 3,$ then, considering the coefficients in front of
$z^{q_{t+1}}$ and $z^{q_{t+1}-1},$ we get
$a_{t+1}(\alpha)=a_{t+1}(0)$ and $b_{t+1}(\alpha)=b_{t+1}(0)$, a
contradiction.

So $q_{t+1}=2.$ If $t>r,$ then for $P_M(z)=z^{t+1}-Mz^{t-1}$ we
have diam$Z(P_M)=2\sqrt{|M|}.$ On the other hand, since
$$L[P_M]=(a_{t+1}(0)z+b_{t+1}(0))^2-M(a_{t-1}(0)z+b_{t-1}(0))$$
it follows that
$$\lim_{M\to\infty}\frac{{\rm diam}Z(L[P_M])}{|M|}=\frac{a_{t-1}(0)}{a_{t+1}(0)}\neq
0$$ which is a contradiction.

Hence $q_{r+1}=2.$ Then
$$(a_{r+1}(\alpha)z+b_{r+1}(\alpha))^2=(a_{r+1}(0)z+b_{r+1}(0))^2
+\alpha(r+1)(a_rz+b_r)+Q_{r+1}(\alpha),$$ where
$$Q_{r+1}(\alpha)=\sum_{j=2}^{r+1}\binom{r+1}{j}\alpha^jL[z^{r+1-j}].$$
Setting $b_{r+1}=b_{r+1}(0),$ it follows that
$$a_{r+1}(\alpha)=a_{r+1}(0)=a_{r+1},$$
$$b_{r+1}(\alpha)=b_{r+1}+\frac{\alpha(r+1)a_r}{2a_{r+1}}\;\;\;{\rm and}$$
$$b^2_{r+1}(\alpha)=b^2_{r+1}+\alpha(r+1)b_r+Q_{r+1}(\alpha).$$
The last two equalities imply that
$\displaystyle\frac{b_{r+1}}{a_{r+1}}=\frac{b_r}{a_r},$ $r\ge 1$,
$$L[z^{r-1}]=\frac{(r+1)a_r^2}{2ra_{r+1}^2}\;\;\;{\rm and}\;\;\;  L[z^{r-2}]=\dots=L[1]=0.$$

Since $L$ is a linear operator, it is enough to consider the case
$L[z^{r-1}]=(r-1)!.$ Set $\displaystyle a=\frac{a_r}{r!},$ $\displaystyle
b=\frac{b_r}{r!}$ and $\mathcal{L}(z)=az+b.$ It follows that $\displaystyle
a_{r+1}=\frac{(r+1)!a^2}{2},$ $\displaystyle b_{r+1}=\frac{(r+1)!b^2}{2}$
and hence $L[P]=a^{1-r}(P\circ\L)^{(r-1)}$ for
$P(z)=z^j,$ $0\le j\le r+1.$

We shall now prove by induction the last equality in general.
Assume that it holds for $j=0,\dots,k$ $(r+1\le k\le n-1).$ We see
as before (by comparing coefficients) that $q_{k+1}=q_k+1$ or
$q_{k+1}=q_k.$ The last case is impossible by considering, as above,
the polynomial $z^{k+1}-Mz^{k-1}.$ So $q_{k+1}=k+2-r\ge 3$ and
$$L[(z+\alpha)^{k+1}]-L[z^{k+1}]
=(a_{k+1}(\alpha)z+b_{k+1}(\alpha))^{q_{k+1}}
-(a_{k+1}(0)z+b_{k+1}(0))^{q_{k+1}}$$
$$=\sum_{j=1}^{k+1}\binom{k+1}{j}\alpha^j L[z^{k+1-j}]=\tilde
P_{k+1}(z+\alpha)-\tilde P_{k+1}(z),$$ where
$$\tilde P_{k+1}(z)=\frac{(k+1)!}{(q_{k+1})!}(az+b)^{q_{k+1}}.$$
In particular, $a_{k+1}(\alpha)=a_{k+1}(0).$
\smallskip

{\it Claim.} If $l\ge 3$ and $\beta\neq 0$ are such that for any
$w\in\mathcal{C},$
$$(w+\beta)^l-w^l=d[(w+\gamma)^l-(w+\delta)^l],$$
then $d=1,\beta=\gamma,\delta=0$ or $d=-1,\beta=\delta,\gamma=0.$
\smallskip

This follows by the fact that if $\lambda_0,\dots,\lambda_l$ are
pairwise different complex number, then
$(w+\lambda_0)^l,\dots,(w+\lambda_l)^l$ is a basis of $\mathcal{P}_l,$

The claim easily implies that $L[z^{k+1}]=\tilde P_{k+1}(z)$, which
completes the induction.

Finally, note that diam$Z(L[P])\le$ diam$Z(P)$ leads to $|a|\ge
d_{n,k}.$

The proof of Proposition 1 is completed.

\end{document}